\documentclass[12pt,a5paper]{article}
\usepackage{url}
\usepackage[colour,hyperref,a5]{ajr}

\title{The prisoners may be in two minds}
\author{A. J. Roberts\thanks{Dept Maths \& Computing, University of 
Southern Queensland, Toowoomba, Queensland 4352, \textsc{Australia}. 
\protect\url{mailto:aroberts@usq.edu.au}}}

\newcommand{\nash}[1]{\textbf{#1}}
\newcommand{\B}[1]{\textit{#1}}
\newcommand{\pa}[1]{{\color{red}#1}}
\renewcommand{\pm}[1]{{\color{blue}#1}}
\newcommand{\qa}[1]{{\color{green}#1}}
\newcommand{\qm}[1]{{\color{magenta}#1}}

\newsavebox{\ajrquotes}
\newenvironment{quotes}[1]{\begin{quote}\savebox{\ajrquotes}{\em #1}}%
{\hfil\penalty50\quad\mbox{}\hfil\usebox{\ajrquotes}\penalty-10000\end{quote}}

\begin{document}

\maketitle

\begin{abstract}
	Recognise that people have many, possibly conflicting, aspects to
	their personality.  We hypothesise that each separate
	characteristic of a personality may be treated as an independent
	player in a non-zero sum many player game.  This idea is applied to
	the two person Prisoners' Dilemma as an introductory example.  We
	assume each prisoner has a ``mercenary'' characteristic as well as
	an ``altruistic'' characteristic, and find that all Nash equilibria of
	the Prisoners' Dilemma has each prisoner in an internal conflict
	between their two characteristics.  The hypothesis that people are
	composed of more than one ``player'' may explain some of the
	anomalies that occur in human experiments exploring game theory.
\end{abstract}

\section{Introduction}

In reviewing the implications of some psychological research in game
playing, Matthew Rabin comments 
\begin{quotes}{Rabin~\cite[p.12]{Rabin98}}
	humans differ from the way they are traditionally described by
	economists \ldots\ it is sometimes misleading to conceptualise
	people as attempting to maximize a coherent, stable and accurately
	perceived [utility]~$U(x)$\,.
\end{quotes}
To include emotions into artificial intelligences, some computer
scientists explore the modelling of human emotions: Smith and Ellworth,
see \cite[\S8.2]{Baillie02} for example, identified a six~dimensional
``affective space'' to capture 15~emotions.\footnote{Smith and
Ellworth's axes for the six dimensional affective space are called:
pleasantness, anticipated effort, certainty, attentional activity,
responsibility and control.  The fifteen categorised emotions are:
happiness, sadness, anger, boredom, challenge, hope, fear, interest,
contempt, disgust, frustration, surprise, pride, shame and guilt.} This
suggests humans are reasonably modelled by six ``appraisal''
dimensions.  Crucially, these appraisal dimensions are independent.
Maybe humans playing games, including the game called economics, act
according to the diverse needs of such internal characteristics of a
complex personality; and not according to maximising a single well
defined utility.

\section{The prisoners' dilemma}

As a first tentative exploration of the idea of multiple game players
within one person, we explore two people, P~and~Q, playing the
Prisoners' Dilemma.  Suppose each person has an altruistic, cooperative
side and a mercenary, uncooperative side to their personality.  These
characteristics are not to be considered as opposites, but as
independent characteristics, independent ``players'', within each
person.  Perhaps view it as each person having two independent
dimensions to their character.\footnote{Of course, we presume that a
more realistic model of personality would have more than two
independent dimensions.}

Thus within person~P we have two independent players, \pa{Pa} and
\pm{Pm}, the pair of altruistic and mercenary characters respectively
within person~P. Similarly \qa{Qa} and \qm{Qm} are independent
characteristics of person~Q{}.  Be careful of the
distinction between people and players: each person is supposed to be a
composite of independent players.  The two people then face a
four-player Prisoner's Dilemma---two players per person.  Each player
has a choice between the two strategies of cooperation~$C$ (remaining
silent) and defecting~$D$ (informing).  We suppose, within any one
person, that if either of the internal players defect, whether the
mercenary or the altruistic, then the person does defect---a chain is
only as strong as its weakest link.  That is, both aspects of a
person's character have to cooperate in order for the person to
actually cooperate.

For simplicity we rate the outcomes on a four point scale of
preferences for each player.  The preferences are not symmetric, but we
suppose symmetry between the pair of mercenary players and between the
pair of altruistic players.  We also assume the preferences of a player
only depend upon the other person's actual action: that is, whether they
actually cooperate or defect.  Thus the preferences only depend upon
the other person's players strategies in the combinations $\{CC\}$ and
$\{CD,DC,DD\}$.  We suppose the following preferences for the two types
of player within each person.
\begin{description}
	\item[Mercenary.]  These players only care about the actual 
	outcome and thus the preferences depend upon his/her persons 
	actual actions, that is, $\{CC\}$ and $\{CD,DC,DD\}$.  As in the 
	classic prisoners dilemma the outcomes in increasing order of 
	preference are:
		\begin{enumerate}
		\item $\{CC\}\times\{CD,DC,DD\}$ when the person cooperates but
		the other defects;
	
		\item  $\{CD,DC,DD\}\times\{CD,DC,DD\}$ when both people defect;
	
		\item  $\{CC\}\times\{CC\}$ when both cooperate;
	
		\item  $\{CD,DC,DD\}\times\{CC\}$ when the person defects and the 
		other cooperates.
	\end{enumerate}
	These mercenary preferences are shown in \pm{blue} and \qm{magenta}
	in Table~\ref{tbl:pref}.
	
\begin{table}[tbp]
	\centering
	\caption{payoff preferences for \pm{P-mercenary},
	\pa{P-altruistic}, \qm{Q-mercenary} and \qa{Q-altruistic}.  The
	best response actions are indicated in \B{italic}, and the four
	Nash equilibria in \nash{\B{bold italics}} are seen in the middle
	of the payoffs.}
\begin{tabular}{|c|c|c|}\hline
	 & $\qm{\mbox{Qm}=C}$ & $\qm{\mbox{Qm}=D}$  \\[1ex] \hline
	$\pm{\mbox{Pm}=C}$ 
	& 
	\begin{tabular}{ccc}
		 & $\qa{\mbox{Qa}=C}$ & $\qa{\mbox{Qa}=D}$  \\
		$\pa{\mbox{Pa}=C}$ & \pm{3}\,\pa{\B3}\,\qm{3}\,\qa{\B3} & \pm{1}\,\pa{\B1}\,\qm{\B4}\,\qa{\B3}  \\
		$\pa{\mbox{Pa}=D}$ & \pm{\B4}\,\pa{\B3}\,\qm{1}\,\qa{\B1} & \nash{\pm{\B2}\,\pa{\B1}\,\qm{\B2}\,\qa{\B1}}  \\
	\end{tabular}
	&   
	\begin{tabular}{cc}
		  $\qa{\mbox{Qa}=C}$ & $\qa{\mbox{Qa}=D}$  \\
		  \pm{1}\,\pa{\B1}\,\qm{\B4}\,\qa{\B4} & \pm{1}\,\pa{\B1}\,\qm{\B4}\,\qa{3}  \\
		  \nash{\pm{\B2}\,\pa{\B1}\,\qm{\B2}\,\qa{\B2}} & \pm{\B2}\,\pa{\B1}\,\qm{\B2}\,\qa{1}  \\
	\end{tabular}
	\\[4ex] \hline
	$\pm{\mbox{Pm}=D}$ 
	&  
	\begin{tabular}{ccc}
		$\pa{\mbox{Pa}=C}$ & \pm{\B4}\,\pa{\B4}\,\qm{1}\,\qa{\B1} & \nash{\pm{\B2}\,\pa{\B2}\,\qm{\B2}\,\qa{\B1}}  \\
		$\pa{\mbox{Pa}=D}$ & \pm{\B4}\,\pa{3}\,\qm{1}\,\qa{\B1} & \pm{\B2}\,\pa{1}\,\qm{\B2}\,\qa{\B1}  \\
	\end{tabular}
	&   
	\begin{tabular}{cc}
		  \nash{\pm{\B2}\,\pa{\B2}\,\qm{\B2}\,\qa{\B2}} & \pm{\B2}\,\pa{\B2}\,\qm{\B2}\,\qa{1}  \\
		  \pm{\B2}\,\pa{1}\,\qm{\B2}\,\qa{\B2} & \pm{\B2}\,\pa{1}\,\qm{\B2}\,\qa{1}  \\
	\end{tabular}
	\\ \hline
\end{tabular}
	\label{tbl:pref}
\end{table}%

	\item[Altruistic.]  These players do care about the outcome as for 
	the mercenary characters, but if they, the altruistic player, 
	defects then guilt lowers their preference of the outcome.  In 
	particular and for simplicity we suppose the preference is 
	decreased to the one lower.  Thus for an altruistic player the 
	preferences in increasing order are:
		\begin{enumerate}
		\item $\{CC,CD,DD\}\times\{CD,DC,DD\}$ when the person 
		cooperates but the other defects, or the altruistic player 
		defects and the other person defects;
	
		\item $\{DC\}\times\{CD,DC,DD\}$ when both people defect but
		only via the mercenary player;
	
		\item  $\{CC,CD,DD\}\times\{CC\}$ when both cooperate, or when the 
		altruistic player defects and the other person cooperates;
	
		\item $\{DC\}\times\{CC\}$ when only the mercenary player 
		defects and the other person cooperates.
	\end{enumerate}
	These altruistic preferences are shown in \pa{red} and \qa{green} in 
	Table~\ref{tbl:pref}.

\end{description}
Inspect the consequent preferences shown in Table~\ref{tbl:pref}.  The
best response~\cite[\S2.8,e.g.]{Osborne02} for each player as as
function of the other three players choices are shown in \B{italic}.  A
Nash equilibrium corresponds to any cell with all four preferences
being such a best response.  See the four Nash equilibria (in the
middle of the table) are obtained from the strategies: 
\begin{displaymath}
    \{\pm C\pa D\qm
C\qa D,\pm C\pa D\qm D\qa C,\pm D\pa C\qm C\qa D,\pm D\pa C\qm D\qa
C\}\,.
\end{displaymath}
Remarkably, the four Nash equilibria correspond to \emph{both
prisoners being in two minds} about what action to take.  Perhaps this
indicates something about the internal stress suffered by a person in a
Prisoner's Dilemma.  The stress comes from the internal conflict
between the different characteristics within each person.

Since the four-player Prisoner's Dilemma only involves preferences the
Nash equilibria are reasonably robust.\footnote{I have not searched for
any mixed Nash equilibria.}

\section{Discussion}


\begin{quotes}{Rabin~\cite[p.16]{Rabin98}}
	Yet pure self-interest is far from a complete description of human
	motivation, and realism suggests that economists should move away
	from the presumption that people are \emph{solely} self-interested.
\end{quotes}
Here we have discussed one model for how generalise the analysis of
human behaviour by positing a complex interplay of motivation and
reward internal to each person.  This is a type of multiple-self
model of human behaviour.

The proposed model of people composed of multiple independent selves
also suggests a rationale for framing effects, described as:
\begin{quotes}{Rabin~\cite[p.36]{Rabin98}}
	two logically equivalent (but not \emph{transparently} equivalent)
	statements of a problem lead decision makers to choose different
	options.
\end{quotes}
The many Nash equilibria that potentially exist in a game by people
with multiple selves make it quite likely that different statements of
a situation will lead to different Nash equilibria being realised.
Recall that the extensive form of a game often favours one Nash
equilibria over another~\cite[\S6.2]{Osborne94}, even when the two Nash
equilibria are equally valid in the strategic form of the game; subgame
perfect equilibria are the rational solutions in an extensive game.
Problem statements which present equivalent information in a different
sequential order lead to the playing of different extensive games even
though the corresponding strategic games are equivalent.  Thus the
framing of a game problem possibly explores independent players within
a person.

A different multiple-self model has been proposed to explain time varying
preferences:
\begin{quotes}{Rabin~\cite[p.39]{Rabin98}}
	a person is modeled as a separate ``agent'' who chooses her current
	behavior to maximize her current long-run preferences, whereas each
	of her future selves, with her own preferences, will choose her
	future behaviour to maximize \emph{her} preferences.
\end{quotes}
The difference is that here we posited multiple selves to co-exist
simultaneously within each person.  Such a multiple-self model could 
explain the richness of human behaviour much better than one simple
utility.

\paragraph{Acknowledgements:} I thank Daniel Burrell, Sam Farrell and
Penny Baillie for their interest in the concept of a multi-player
person.

\bibliographystyle{plain}
\bibliography{ajr,new,bib}

\end{document}